\documentclass[reqno,11pt]{article}

\usepackage{mathtools}
\usepackage{amsmath}
\usepackage{amsfonts}
\usepackage{amssymb}
\usepackage{graphicx}
\usepackage{epstopdf}

\usepackage[inner=2cm, outer=2cm, bottom=2cm, top=2cm]{geometry}

\newcommand*{\eqb}{\begin{equation}}
\newcommand*{\eqe}{\end{equation}}
\newcommand*{\al}{\alpha}

\newcommand{\abs}[1]{\left| #1 \right|}

\newcommand{\norm}[1]{\left\| #1 \right\|}

\newcommand{\scalar}[1]{\left\langle #1 \right\rangle}

\newcommand{\R}{\mathbb{R}}

\newcommand{\pr}{\mathbf{P}}
\newcommand{\ex}{\mathbf{E}}
\newcommand{\N}{\mathbb{N}}

\newcommand{\eps}{\varepsilon}

\newcommand{\bol}[1]{\mathbf{#1}}

\begin{document}

\title{Method of calculating densities for isotropic L\'evy walks}

\author{Marcin Magdziarz $^{1,*}$, Tomasz Zorawik $^1$\\ \ \\
\footnotesize
$^1$ Hugo Steinhaus Center, \\
Faculty of Pure and Applied Mathematics, \\
Wroclaw University of Science and Technology, \\
Wyspianskiego 27, 50-370 Wroclaw, Poland. \\[4pt]
$^{*}$\textit{Author for correspondence (marcin.magdziarz@pwr.edu.pl)} }

\date{}

\maketitle

\abstract{
We provide explicit formulas for asymptotic densities of $d$-dimensional isotropic L\'evy walks,  when $d>1$. The densities of multidimensional undershooting and overshooting L\'evy walks are presented as well. Interestingly, when the number of dimensions is odd the densities of all these L\'evy walks are given by elementary functions. When $d$ is even, we can express the densities as fractional derivatives of  hypergeometric functions, which makes an efficient numerical evaluation  possible.
}


\section{Introduction}
L\'evy walks are one of the most important tools for modeling anomalous stochastic transport phenomena. The period of intensive research on this topic started with the pioneering papers \cite{shlesinger klafter} by Shlesinger et. al. and
\cite{klafter blumen} by Klafter et al. Since then L\'evy walks found many applications in  different areas of physics and biology. The list of real-life phenomena and complex systems where L\'evy walks are used includes (but is not limited to) migration of swarming bacteria \cite{bacteria}, blinking nanocrystals \cite{blinking crystals}, light transport in optical materials \cite{light}, fluid flow in a rotating annulus \cite{annul}, foraging patterns of animals \cite{for1,for2,for3}, human travel \cite{human 1, human 2} and epidemic spreading \cite{epid1,epid2}. For more background information about L\'evy walks and their applications the interested reader is referred to the recent review \cite{levy walks}.

L\'evy walks exhibit two important features. The first one is a power-law jump distribution and the second one is a finiteness of all moments. This is obtained by introducing a dependence between waiting times and lengths of the jumps - we require that the length of the jump  be equal to the preceding waiting time in the underlying continuous-time random walk (CTRW) scenario.  The continuity of the trajectories is obtained by the linear interpolation of the corresponding CTRW. As a result the velocity $v$ of the particle performing L\'evy walk is constant. The appearance of this two features together - the power-law jump distribution and the finiteness of all moments -  stays in contrast with a different very popular model for anomalous transport, namely L\'evy flights \cite{metzler2000, klafter book, JW}. For L\'evy flights the power-law jump distribution implies that the mean square displacement is infinite.
For other correlated fractional diffusion models we refer the reader to \cite{proc1,proc2,proc3,proc4}

In spite of the long history and popularity of L\'evy walks, their multidimensional probability density functions (PDFs) were not known.
In this paper we fill this gap and derive explicit PDFs of $d$-dimensional  L\'evy walks, where $d>1$. When $d$ is an odd number $d=2n+3$ ($n\in\N$), the PDF is expressed by elementary functions. This fact is somehow unexpected since the limit process for L\'evy walks is given as a composition of certain $\al$-stable processes, and it is a known fact that the stable densities are expressed by elementary functions only for special values of $\al$ \cite{JW}. Here the PDFs of L\'evy walks will be expressed in terms of elementary functions for all $\al \in (0,1)$. In the case when the dimension is even $d=2n+2$ ($n\in\N$) the PDF can also be computed, but the formula involves hypergeometric functions and the Riemann-Liouville right-side fractional derivative, which can be efficiently evaluated numerically.  We also provide explicit formulas for the densities of other coupled CTRWs - the so-called undershooting and overshooting L\'evy walks (also known as wait-first and jump-first L\'evy walks \cite{levy walks}). Similarly, these densities are given by elementary functions for odd dimensions $d$.
Moreover these PDFs solve certain differential equations \cite{MSSZ} with the fractional material derivative \cite{material derivative 1, material derivative 2}.

The densities of 1-dimensional ballistic L\'evy walks have been found by Froemberg et al. in \cite{asymptotic densities}, see also \cite{mz_fractional} for other approach to this problem.
Recently, in \cite{tzmm} the PDFs of  $2$ and $3$-dimensional ballistic L\'evy walks were derived.
Comparison of  three different models of L\'evy walks in two dimensions can be found in \cite{Barkai2}.
In this paper we calculate the densities of $d$-dimensional L\'evy walks when $d>1$ is arbitrary. We also apply our method to the overshooting and undershooting L\'evy walks.
The main idea behind this method is to take advantage of the rotational invariance of the L\'evy walks and connect the multi-dimensional PDF with a proper one-dimensional distribution using methods from \cite{uchaikin zolotarev book}, then apply the formula of Godr\`{e}che and Luck from \cite{godreche luck} to invert the Fourier-Laplace transform.

\section{L\'evy walks and their limits}
In this section we recall the definition of $d$-dimensional standard, undershooting and overshooting L\'evy walks (see \cite{levy walks}). We also recall their limit processes.
\subsection{Definition of L\'evy walks}
Let $T_i$ ($i \in \N$) be a sequence of waiting times. We assume that $T_i$ are independent, identically distributed (IID) positive random variables with power-law distribution $P(T_i>t)\propto t^{-1-\al}$, $\al\in(0,1)$. Denote by $N(t)=\max\{k\geq 0:\sum_{i=1}^k T_i \leq t\}$ the corresponding process counting the number of jumps up to time $t$.
Next, let us define the sequence of consecutive jumps
\[
\bol{X}_i=v T_i \bol{V}_i,\;\;\; i=1,2,... \;.
\]
Here, $\{\bol{V}_i\}$ is a sequence of IID random unit vectors distributed uniformly
on the $d-1$ - dimensional hypersphere $\mathbb{S}^{d-1}$. Each vector $\bol{V}_i$ governs the direction of
$i$-th jump. The constant $v$ is the velocity, for simplicity it is assumed that $v=1$. In the Introduction we mentioned that for L\'evy walk the length of each jump should be equal to the corresponding waiting time. It is clear from the above equation that this condition is satisfied: $\norm{\bol{X}_i}=T_i$, where $\norm{}$ denotes the Euclidean norm in $\R^d$.
Now, the undershooting L\'evy walk (or wait-first L\'evy walk) is defined as
\eqb\label{under}
\bol{L}_{ULW}(t)=\sum_{i=1}^{N_t}\bol{X}_i.
\eqe
This is a CTRW, so the trajectories
 are piecewise constant and have jumps.  The overshooting L\'evy walk (or jump-first L\'evy walk) has the following definition
\eqb
\label{over}
\bol{L}_{OLW}(t)=\sum_{i=1}^{N_t+1}\bol{X}_i.
\eqe
This is also a CTRW so the trajectories of this process are also not continuous.
However, applying simple linear interpolation on the trajectories of $\bol{L}_{ULW}(t)$ and $\bol{L}_{OLW}(t)$, we arrive at the final definition of the standard L\'evy walk $\bol{L}(t)$:
\eqb
\begin{split}
\label{definition lw}
\bol{L}(t)&= \bol{L}_{ULW}+\left(t-\sum_{i=1}^{N(t)}T_i\right)(\bol{L}_{OLW}-\bol{L}_{ULW})\\
&=\sum_{i=1}^{N(t)}T_i\bol{V}_i+\left(t-\sum_{i=1}^{N(t)}T_i\right)\bol{V}_{N(t)+1},
\end{split}
\eqe
The trajectories of $\bol{L}(t)$ are continuous and piecewise linear, which means that the walker moves with constant velocity $v$. The spatio-temporal coupling ensures that all moments, in particular a mean square displacement, are finite. This follows from that fact that $\norm{\bol{L}(t)}\leq t$. Similarly, all the  moments of $\bol{L}_{ULW}(t)$ are also finite. However the overshooting L\'evy walk does not have this property, that is $\ex \norm{\bol{L}_{ULW}(t)}=\infty$ for all $t>0$. Since the processes $\bol{L}_{ULW}(t)$, $\bol{L}_{OLW}(t)$  and $\bol{L}(t)$ are $d$-dimensional, in the following we will use the notation $\bol{L}_{ULW}(t)=(L_{ULW_1}(t),...,L_{ULW_d}(t))$, $\bol{L}_{OLW}(t)=(L_{OLW_1}(t),...,L_{OLW_d}(t))$ and $\bol{L}(t)=(L_1(t),...,L_d(t))$ for their coordinates.

\subsection{Limit processes}
Let $\bol{M}_\al(t)$ be the $d$-dimensional rotationally invarinat $\al$-stable process \cite{JW} with the Fourier transform
\eqb
\ex(\exp\{-ik\bol{M}_\al(t)\})=\exp\{-t\norm{k}^\al\}.
\eqe
Moreover, let $U_\al(t)$ be the $\al$-stable subordinator, i.e. one-dimensional strictly increasing $\al$-stable L\'evy process.
The inverse $\al$-stable subordinator is defined as the first passage time of $U_\al (t)$, that is
 $S_\al(t)=\inf\{\tau : U_\al(\tau)>t\}$.  We assume that $\bol{M}_\al(t)$ and $U_\al(t)$ have jumps of equal length $\norm{\Delta M_\al(t)}=\Delta U_\al(t)$ for all $t>0$ almost surely. This assumption can be formally expressed in terms of L\'evy measure of the process $(\bol{M}_\al(t), U_\al(t))$ (see \cite{MT}, \cite{TZM}). Then we have the following convergence of L\'evy walks $\bol{L}(t)$, $\bol{L}_{ULW}(t)$ and $\bol{L}_{OLW}$ in Skorokhod space (which also implies convergence of all finite-dimensional distriubtions).

I. Undershooting L\'evy walk limit \cite{MT}, \cite{TZM}:
\eqb
\frac{\bol{L}_{ULW}(nt)}{n}\xrightarrow{n\rightarrow \infty} \bol{Y}(t),
\eqe
where $\bol{Y}(t)$ is a right continuous version of $\bol{M}_\al^{-}(S_\al^{-1}(t))$. Here $\bol{M}_\al^{-}(t)$ denotes the left limit process $\bol{M}_\al^{-}(t)=\lim_{s\rightarrow t^-} \bol{M}_\al(s)$.

II. Overshooting L\'evy walk limit \cite{MT}:
\eqb
\frac{\bol{L}_{OLW}(nt)}{n}\xrightarrow{n\rightarrow \infty} \bol{Z}(t),
\eqe
where $\bol{Z}(t)=\bol{M}_\al^{-}(S_\al^{-1}(t))$

III. Standard L\'evy walk limit \cite{MSSZ}:
\eqb
\label{z_representation}
\frac{\bol{L}(nt)}{n}\xrightarrow{n\rightarrow\infty} \bol{X}(t)=\left\{\begin{array}{lr}
\bol{Y}(t) &\mbox{if} \quad t\in\mathcal{R}(U_\al) \\
\bol{Y}(t)+\frac{t-G(t)}{H(t)-G(t)}(\bol{Z}(t)-\bol{Y}(t) &\mbox{if} \quad t\notin\mathcal{R}(U_\al),
\end{array}\right.
\eqe
where
\[
\mathcal{R}(U_\al)=\{U_\al(t)\;:\;t\geq0\}\subset[0,\infty).
\]
Moreover
\[
G(t)=U^-_\al(S_\al(t+))
\]
is the moment of the previous jump of $\bol{Y}(t)$ before time $t$ and
\[
H(t)=U_\al(S_\al(t+))
\]
is the moment of the next jump of $\bol{Z}(t)$ after time $t$. Here $S_\al(t+)=\lim_{s\rightarrow t^+}S_\al(s)$.

The trajectories of the limit process $\bol{X}(t)$ are continuous whereas the trajectories of  $\bol{Y}(t)$ and $\bol{Z}(t)$ are discontinuous.

The jump directions $\bol{V}_i$ are uniformly distributed
on the $d-1$-hypersphere $\mathbb{S}^{d-1}$ which implies that the distributions of L\'evy walks $\bol{L}_{ULW}(t)$, $\bol{L}_{OLW}(t)$ and $\bol{L}(t)$  are rotationally invariant - each direction of the motion is equally possible. Therefore to determine the PDF of $\bol{L}(t)$ it is enough to determine the PDF of the radius $\norm{\bol{L}(t)}=\sqrt{L_1^2{(t)}+...+L_d^2{(t)}}$. Similarly to calculate the PDF of undershooting or overshooting L\'evy walk it suffices to determine the PDF of their radiuses.  The same is true for the limit processes of L\'evy walks $\bol{X}(t)$, $\bol{Y}(t)$, $\bol{Z}(t)$ - they are also isotropic and hence it is enough to find the PDFs of their radiuses.

In the next sections we will find the PDFs of $\bol{X}(t)$, $\bol{Y}(t)$, $\bol{Z}(t)$.
We will perform two steps. The first one is to determine the PDF of $X_1(t)$ -  the projection of $\bol{X}(t)$ on the first axis. The second one is to relate the PDFs of  $\norm{\bol{X}(t)}$ and $X_1(t)$. The same ideas are used for undershooting and overshooting L\'evy walk limits.

\section{PDF of standard L\'evy Walk}
Let $H(\bol{x},t)$, $\bol{x}=(x_1,x_2,...,x_d)\in \R^d$.  be the PDF of $\bol{X}(t)=(X_1(t),X_2(t),...,X_d(t))$. From \cite{MSSZ} we know that the  Fourier-Laplace transform of $H(\bol{x},t)$ is given by
\eqb \nonumber
H(\bol{k},s)=\frac{1}{s}g\left(\frac{i\bol{k}}{s}\right)=\frac{1}{s}\frac{\int_{\bol{S}^{d-1}}\left(1-\scalar{\frac{i\bol{k}}{s}, \bol{u}}\right)^{\alpha-1} K(d\bol{u})}{\int_{\bol{S}^{d-1}}\left(1-\scalar{\frac{i\bol{k}}{s}, \bol{u}}\right)^\alpha K(d\bol{u})},
\eqe
where $\bol{k}=(k_1,k_2,...,k_d)\in\R^d$ is the Fourier space variable, $s$ is the Laplace space variable, $K(d\bol{u})$ is the uniform distribution on a hypersphere $\bol{S}^{d-1}$ and $ \scalar{\;,\;}$ denotes the standard inner product in $\R^d$. Now we can notice that
the Fourier-Laplace transform of the PDF of $X_1(t)$ is equal to
$H_1(k_1,s)=H((k_1,0,...,0),s)$. The marginal distribution of the uniform distribution on a hypersphere $K_1(d\bol{u})$ has the density
\eqb
K(du_1)=c_d\left(1-u_1^2\right)^{(d-3)/2} du_1,
\eqe
where $c_d=\frac{1}{\sqrt{\pi}}\frac{\Gamma(d/2)}{\Gamma((d-1)/2)}$.
Hence
$
H_1(k_1,s)=\frac{1}{s}g_1\left(\frac{ik_1}{s}\right),
$
where
\eqb
\label{g1_3d}
g_1(\xi)=\frac{\int_{-1}^1(1-\xi u)^{\alpha-1}\left(1-u^2\right)^{(d-3)/2}du}{\int_{-1}^1(1-\xi u)^\alpha \left(1-u^2\right)^{(d-3)/2}du}.
\eqe
We also have the ballistic scaling
\eqb
H_1(x,t)=\frac{1}{t}\Phi_1\left(\frac{x}{t}\right),
\eqe
where $\Phi_1(x)=H_1(x,1)$ is the PDF of $X_1(1)$ and $x\in\R$.

\subsection{Odd number of dimensions - $d=2n+3$}
We start with the case $d=2n+3$. Now we are in position to apply the method of inverting the Fourier-Laplace transform used in \cite{asymptotic densities,godreche luck}, which will allow us to get an explicit formula for $\Phi_1$.  This technique is based on a special representation of the function $g_1(\xi)$ and the Sokhotsky-Weierstrass theorem. This technique works only for one dimensional ballistic processes. It was introduced by Godr\`che and Luck in \cite{godreche luck}  for inverting double Laplace transform and then generalized by Froemberg et al. in \cite{asymptotic densities}  for the Fourier-Laplace transform. After its application we obtain (see Appendix A)
\eqb
\label{formula_phi}
\Phi_1(x)=-\frac{1}{\pi}\lim_{\eps \rightarrow 0} \operatorname{Im}\left[\frac{1}{x+i\eps}g_1\left(-\frac{1}{x+i\eps}\right)\right].
\eqe
After calculating the integrals and taking the imaginary part we obtain
\eqb
\begin{split}
\label{odd_d_phir}
&\Phi_1(x)=\frac{\sin(\pi\al)}{\pi x}\mathcal{B}_1\cdot\bigg\{ (-1)^{j_1+m_2}(1/x+1)^{m_1+j_1+\al}(1/x-1)^{m_2+j_2+\al}(1/x^{2}-1)^{2n-m_1-m_2}\\
&\times\left(\frac{1}{m_1+j_1+\al}\frac{1}{m_2+j_2+\al+1}(1/x-1)+ \frac{1}{m_2+j_2+\al}\frac{1}{m_1+j_1+\al+1}(1/x+1)\right)\bigg\}\\
&/\mathcal{B}_2\cdot\bigg\{ \left((-1)^{j_1+j_2}(1+1/x)^{m_1+m_2+j_1+j_2+2\al+2}+(-1)^{m_1+m_2}(1/x-1)^{m_1+m_2+j_1+j_2+2\al+2} \right.\\
 &\left.+2\cos(\pi\al)(-1)^{j_1+m_2}(1/x+1)^{m_1+j_1+\al+1}(1/x-1)^{m_2+j_2+\al+1}\right)(1/x^2-1)^{2n-m_1-m_2}\bigg\},
\end{split}
\eqe
for $x\in(0,1)$, $\Phi_1(-x)=\Phi_1(x)$ for $x\in(-1,0)$ and $\Phi_1(x)=0$ otherwise. The symbols $\mathcal{B}_1$ and $\mathcal{B}_2$ are given by
\eqb
\mathcal{B}_1=\sum_{m_1,m_2=0}^n\sum_{j_1=0}^{m_1}\sum_{j_2=0}^{m_2}{n \choose m_1}{m_1 \choose j_1}{n \choose m_2}{m_2 \choose j_2}2^{m_1+m_2-j_1-j_2},
\eqe
\eqb
\mathcal{B}_2=\sum_{m_1,m_2=0}^n\sum_{j_1=0}^{m_1}\sum_{j_2=0}^{m_2}{n \choose m_1}{m_1 \choose j_1}{n \choose m_2}{m_2 \choose j_2}2^{m_1+m_2-j_1-j_2}\frac{1}{(m_1+j_1+\alpha+1)(m_2+j_2+\alpha+1)}.
\eqe
Alternatively, we can express the density $\Phi_1$ in terms of the hypergeometric functions $_2F_1(a,b;c;x)$  \cite{hiper}. Then
\eqb
\label{odd_d_phir_hyper}
\Phi_1(x)=-\frac{1}{\pi \abs{x}}\operatorname{Im}\frac{_{2}F_1((1 - \alpha)/2, 1 - \alpha/2; 3/2+n; \frac{1}{x^2})}{_2F_1(-\alpha/2, (1- \alpha)/2; 3/2+n; \frac{1}{x^2})}.
\eqe

We will now find the relationship between the PDF $H_R(r,t)$ of the radius $\norm{\bol{X}(t)}$ and the already calculated PDF of $X_1(t)$. As we already noted, the process $\bol{X}(t)$ is isotropic. Hence to find $H(\bol{x},t)$ we only have to determine $H_R(r,t)$. From the property $H(\bol{k},s)=\frac{1}{s}g\left(\frac{i\bol{k}}{s}\right)$ we get the scaling  $H_R(r,t)=\frac{1}{t}\Phi_R\left(\frac{r}{t}\right)$, where $\Phi_R(r)$ is the PDF of
$\norm{\bol{X(1)}}$.
The factorization of $\bol{X}(1)$ into radial
and directional parts gives us
\eqb
\bol{X}(1)\stackrel{d} =\norm{\bol{X(1)}}\bol{V},
\eqe
where $\bol{V}$ is a random vector uniformly distributed on a $d-1$-hypersphere, independent of $\norm{\bol{X(1)}}$ and "$\stackrel{d}=$" denotes the equality of distribution. We infer that
\eqb
\label{2d relation}
\pr (|X_1(1)|\leq x)=2c_n\int_0^1\pr\left(\norm{\bol{X}(1)}\leq\frac{x}{u}\right)\left(1-u^2\right)^{n}du
\eqe
where $c_n=\frac{1}{\sqrt{\pi}}\frac{\Gamma(n+3/2)}{\Gamma(n+1)}$ for $x\geq0$. After differentiation
\eqb
\Phi_1(x)=c_n\int_0^1(1-u^2)^{n}\frac{1}{u}\Phi_R\left(\frac{x}{u}\right)du.
\eqe
Hence
\eqb
\Phi_1(\sqrt{x})=\frac{c_n\Gamma(n+1)}{2} I_-^{n+1}\{\Phi_R(\sqrt{u})u^{-n-1}\}(x),
\eqe
where $I_-^{n+1}$ denotes the right-side Riemann-Liouville integral $I_-^{\beta}$ of the order $\beta=n+1$ \cite{frac}:
\eqb
\label{rl_integral}
I_-^{\beta}\{f(x)\}(y)=\frac{1}{\Gamma(\beta)}\int_x^\infty\frac{f(t)dt}{(t-x)^{1-\beta}}.
\eqe
We can now recover $\Phi_R$:
\eqb
\label{odd relation}
\Phi_R(\sqrt{r})=\frac{2\sqrt{\pi}}{\Gamma (n+3/2)}r^{n+1}(-1)^{n+1}\frac{d^{n+1}}{dr^{n+1}}\Phi_1(\sqrt{r}).
\eqe
if $r\in(0,1)$ and $\Phi_R(r)=0$ otherwise. Since $\Phi_1$ is expressed via elementary functions, so is $\Phi_R$. The densities $\Phi_R$ obtained with this formula for different values of $n$ are plotted in Fig. 2.
For example, in the case $d=3$ which corresponds to $n=0$ we recover the result for a radius from \cite{tzmm}:
\begin{eqnarray}
\label{3d radius}
&&\Phi_R(r)=\frac{8}{\pi}\frac{\al+1}{\al}\sin(\pi\al)r(1 - r^2)^{\al-1} \times\nonumber\\
&&\frac{(1 + r)^{2+2\al}(1+\al - r) - (1-r)^{2+2\al}(1+\al + r)-2r(1-r^2)^{1+\al}\cos(\pi\al)}
{\left((1 + r)^{2+2\al} + (1 - r)^{2+2\al}+2(1-r^2)^{1+\al}\cos(\pi\al)\right)^2}.
\label{3d result}
\end{eqnarray}
In Cartesian coordinates $H(\bol{x},t)$ can be calculated as
\eqb
\label{odd cartesian}
H(\bol{x},t)=\frac{\Gamma(n+3/2)}{2\pi^{n+3/2} t\norm{\bol{x}}^{2n+2}}\Phi_R\left(\frac{\norm{\bol{x}}}{t}\right),
\eqe
with $\norm{\bol{x}}=\sqrt{x_1^2+x_2^2+...+x_{d}}$. We underline that the final result is given by elementary function for any $n$. However, due to the multiple differentiation $(n+1$ times) of the quotient, it is long and hence difficult to write it in one close equation.
The density $H(\bol{x},t)$ of $3$-dimensional standard L\'evy walk is plotted for $t=1$ in Fig. 1.
\begin{figure}
\begin{center}
\includegraphics[width=10cm]{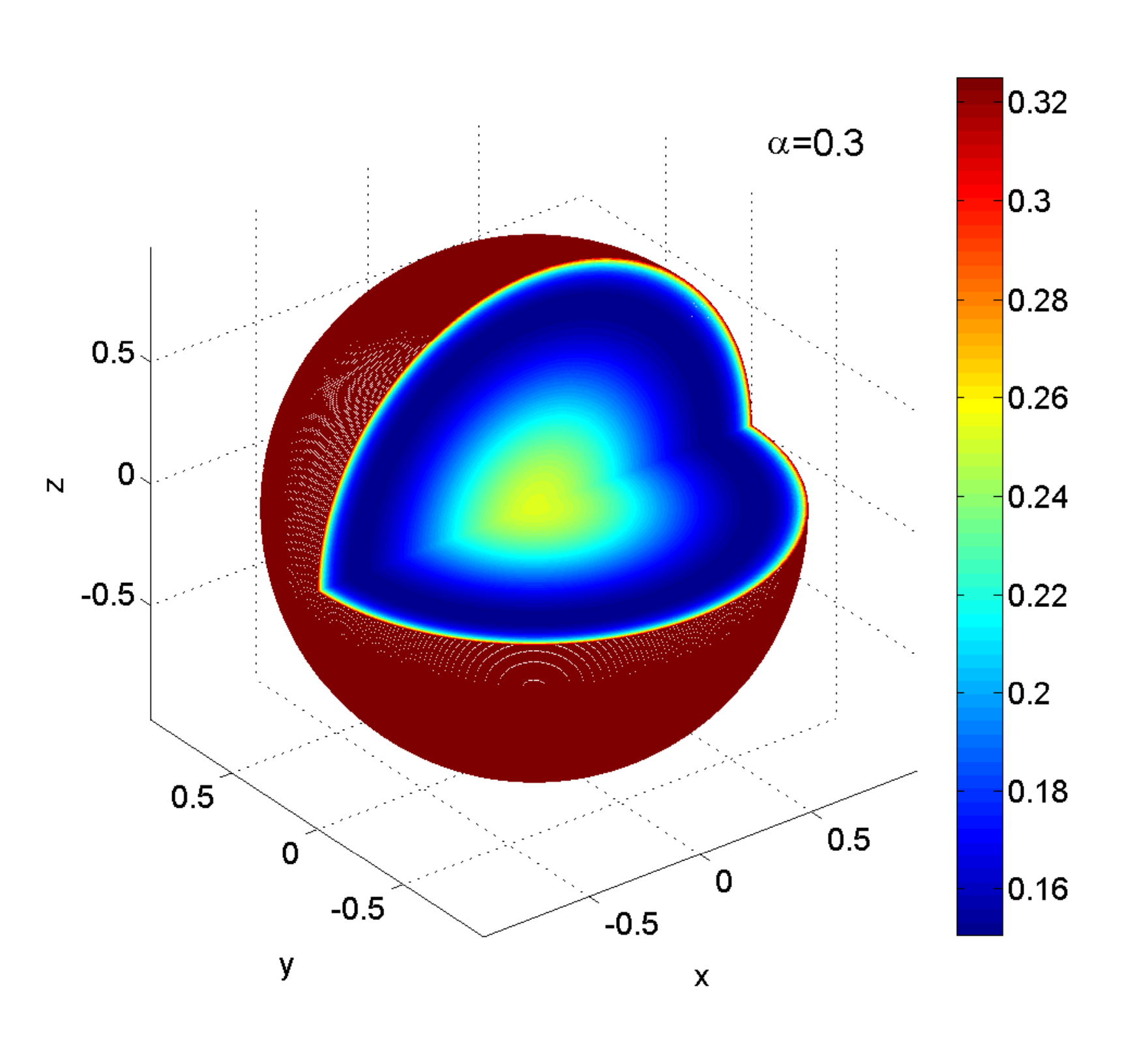}
\label{sphere3d}
\caption{Density $H(\bol{x},t)$ of $3$-dimensional isotropic standard L\'evy walk for $t=1$ and $\alpha=0.3$ obtained from Eqs. \ref{3d radius} and \ref{odd cartesian}. }
\end{center}
\end{figure}
\subsection{Even number of dimensions - $d=2n+2$}
The calculations are getting more complicated when $d$ is even, however the idea remains the same. First we find the PDF $\Phi_1(x)$ of the first-coordinate process $X_1(1)$. The function $g_1(\xi)$ (Eq. \eqref{g1_3d}) now has the form
\eqb
\label{g1 even}
g_1(\xi)=\frac{\int_{-1}^1(1-\xi u)^{\alpha-1}\left(1-u^2\right)^{n-1/2}du}{\int_{-1}^1(1-\xi u)^\alpha \left(1-u^2\right)^{n-1/2}du}.
\eqe
 Using the hypergeometric functions the above equation can be rewritten:
\eqb
g_1(\xi)=\frac{_{2}F_1((1 - \alpha)/2, 1 - \alpha/2; 1+n; \xi^2)}{_2F_1(-\alpha/2, (1- \alpha)/2; 1+n; \xi^2)}.
\eqe
Once again we apply the method of Godr\`{e}che and Luck to get
\eqb
\label{even_d_phir}
\Phi_1(x)=-\frac{1}{\pi \abs{x}}\operatorname{Im}\frac{_{2}F_1((1 - \alpha)/2, 1 - \alpha/2; 1+n; \frac{1}{x^2})}{_2F_1(-\alpha/2, (1- \alpha)/2; 1+n; \frac{1}{x^2})},\quad x\in (0,1).
\eqe
The relation between $\Phi_R$ and $\Phi_1$ now has the form
\eqb\pr (|X_1(1)|\leq x)=2c_n\int_0^1\pr\left(\norm{\bol{X}(1)}\leq\frac{x}{u}\right)\left(1-u^2\right)^{n-1/2}du,
\eqe
where $c_n=\frac{1}{\sqrt{\pi}}\frac{\Gamma(n+1)}{\Gamma(n+1/2)}$ for $x\geq0$. Following the reasoning we used for the odd number of dimensions we obtain
\eqb
\Phi_1(\sqrt{x})=\frac{c_n\Gamma(n+1/2)}{2} I_-^{n+1/2}\{\Phi_R(\sqrt{u})u^{-n-1}\}(x)
\eqe
and hence
\eqb
\label{even relation}
\Phi_R(\sqrt{r})=\frac{2\sqrt{\pi}}{\Gamma (n+1)}r^{n+1/2}D_-^{n+1/2}\{\Phi_1(\sqrt{t})\}(r)
\eqe
if $r\in(0,1)$ and $\Phi_R(r)=0$ otherwise.
Here $D_-^{n+1/2}$ is the right-side Riemann-Liouville fractional derivative of order $n+1/2$ \cite{frac}:
\eqb
\label{rl_derivative}
D_-^{n+1/2}\{f(x)\}(y)=\left(-\frac{d}{d y}\right)^n\frac{1}{\pi^{1/2}}\int_y^\infty \frac{f(x)}{(x-y)^{1/2}}dx.
\eqe
Fig. 2 presents the densities $\Phi_R$ for different number of dimensions $d$. For even $d$ it is necessary to numerically calculate the fractional derivative. This was done in Matlab using algorithm from \cite{matlab}.
In Cartesian coordinates $H(\bol{x},t)$ can be calculated as
\eqb
H(\bol{x},t)=\frac{\Gamma(n+1)}{2\pi^{n+1} t\norm{\bol{x}}^{2n+1}}\Phi_R\left(\frac{\norm{\bol{x}}}{t}\right).
\eqe
\begin{figure}
\begin{center}
\includegraphics[width=10cm]{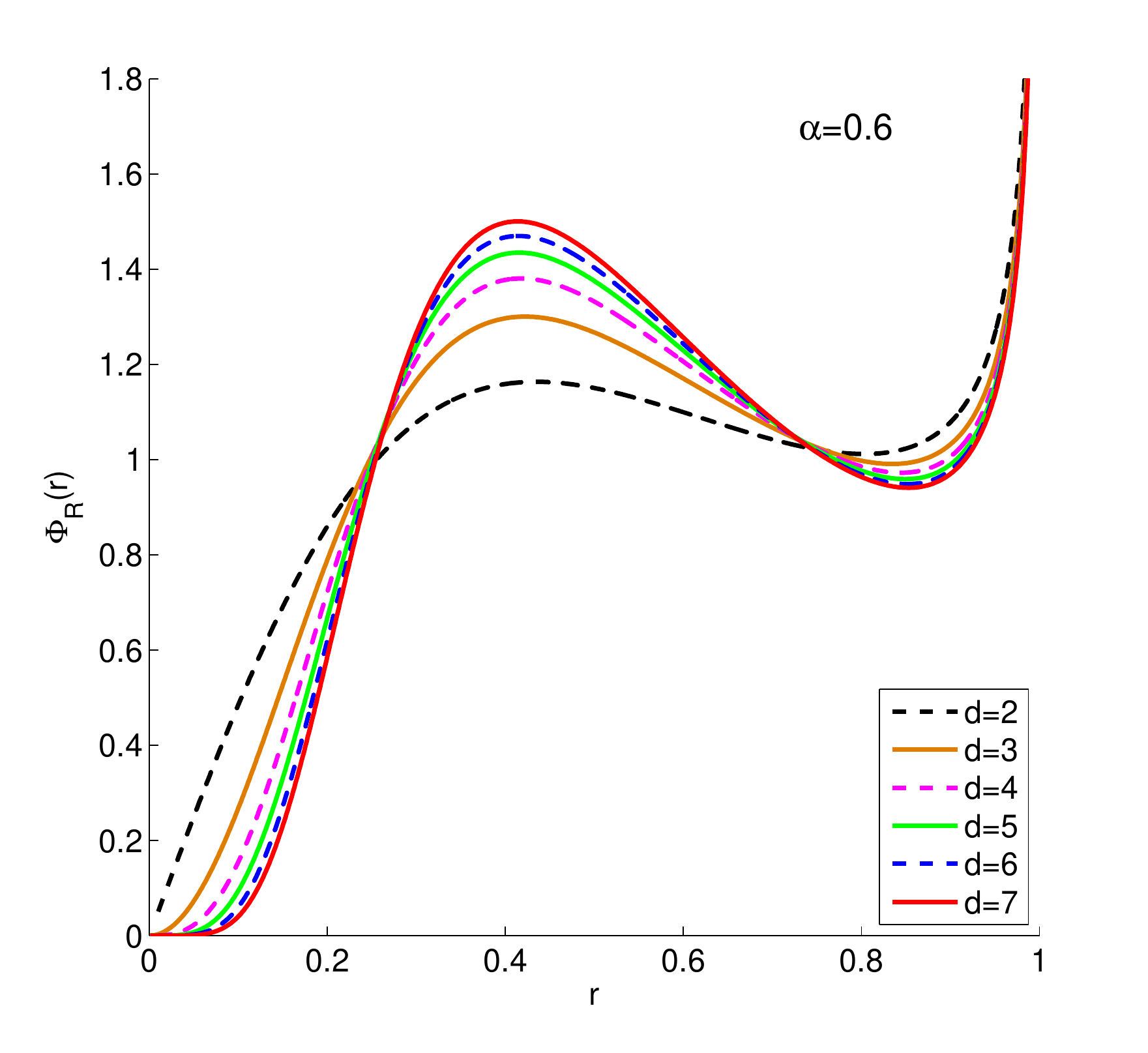}
\label{mcont}
\caption{Densities $\Phi_R(r)$ of radii of $d$-dimensional isotropic standard L\'evy walks for $\alpha=0.6$. Odd number of dimensions - Eq. \ref{odd relation} (elementary functions), even number -  Eq. \ref{even relation} and numerical scheme \cite{matlab}.}
\end{center}
\end{figure}

%
%
\section{Undershooting and overshooting L\'evy walks}
The methods which we used to get the explicit densities of standard L\'evy walks can be also applied for limits of under- and overshooting L\'evy walks.
Recall the notation $\bol{Y}(t)=(Y_1(t),Y_2(t),...,Y_d(t))$ for undershooting and $\bol{Z}(t)=(Z_1(t), Z_2(t),..., Z_d(t)$ for overshooting L\'evy walk.

The relation between the PDF $\Phi_1^Y$ of $Y_1(1)$ and the PDF $\Phi_R^Y$ of $\norm{\bol{Y}(1)}$ remains the same as for standard L\'evy walks. Similarly for the PDF $\Phi_1^Z$ of $Z_1(1)$ and the PDF $\Phi_R^Z$ of $\norm{\bol{Z}(1)}$. We only have to calculate $\Phi_1^Y$ and $\Phi_1^Z$.
\subsection{Undershooting - odd number of dimensions $d=2n+3$}
The PDF of $\bol{Y}(t)$ in Fourier-Laplace space is given by (see \cite{MT})
\eqb \nonumber
H(\bol{k},s)=\frac{1}{s}g\left(\frac{i\bol{k}}{s}\right)=\frac{1}{s}\frac{1}{\int_{\bol{S}^{d-1}}\left(1-\scalar{\frac{i\bol{k}}{s}, \bol{u}}\right)^\alpha K(d\bol{u})},
\eqe
where $\bol{k}=(k_1,k_2,...,k_d)\in\R^d$ is the Fourier space variable, $s$ is the Laplace space variable. The inversion formula yields
\eqb
\begin{split}
&\Phi^Y_1(x)=\frac{\sin(\pi\al)}{\sqrt{\pi}x^{2n+2}}\frac{\Gamma(n+1)}{\Gamma(n+3/2)}\mathcal{B}_1\cdot\bigg\{(-1)^m(1/x-1)^{m+j+\al+1}(1/x^{2}-1)^{n-m}\bigg\}\\
&/\mathcal{B}_2 \cdot\bigg\{\left((-1)^{j_1+j_2}(1+1/x)^{m_1+m_2+j_1+j_2+2\al+2}+(-1)^{m_1+m_2}(1/x-1)^{m_1+m_2+j_1+j_2+2\al+2} \right.\\
 &\left.+2\cos(\pi\al)(-1)^{j_1+m_2}(1/x+1)^{m_1+j_1+\al+1}(1/x-1)^{m_2+j_2+\al+1}\right)(1/x^{2}-1)^{2n-m_1-m_2}\bigg\},
\end{split}
\eqe
for $x\in(0,1)$, $\Phi^Y_1(-x)=\Phi^Y_1(x)$ for $x\in(-1,0)$ and $\Phi^Y_1(x)=0$ otherwise. The symbols $\mathcal{B}_1$ and $\mathcal{B}_2$ equal here
\eqb
\mathcal{B}_1=\sum_{m}^n\sum_{j=0}^{m}{n \choose m}{m \choose j}2^{m-j}\frac{1}{m+j+\alpha+1},
\eqe
\eqb
\mathcal{B}_2=\sum_{m_1,m_2=0}^n\sum_{j_1=0}^{m_1}\sum_{j_2=0}^{m_2}{n \choose m_1}{m_1 \choose j_1}{n \choose m_2}{m_2 \choose j_2}2^{m_1+m_2-j_1-j_2}\frac{1}{(m_1+j_1+\alpha+1)(m_2+j_2+\alpha+1)}.
\eqe
One can also use special functions
\eqb
\begin{split}
\label{odd_under_d_phir_hyper}
\Phi^Y_1(x)=-\frac{1}{\pi \abs{x}}\operatorname{Im}\frac{1}{_2F_1(-\alpha/2, (1- \alpha)/2; 3+n/2; \frac{1}{x^2})}.
\end{split}
\eqe
The PDF $\Phi^Y_R$ is calculated from the equation
\eqb
\Phi_R^Y(\sqrt{r})=\frac{2\sqrt{\pi}}{\Gamma (n+3/2)}r^{n+1}(-1)^{n+1}\frac{d^{n+1}}{dr^{n+1}}\Phi_1^Y(\sqrt{r})
\eqe
if $r\in(0,1)$ and $\Phi^Y_R(r)=0$ otherwise.
In the special case $d=3$ this gives us
\eqb
\begin{split}
&\Phi^Y_R(r)=\frac{4(\alpha+1)}{\pi}\sin(\pi\alpha)(1 - r)^{\alpha}r^{-1 + \alpha}\times\\
&\frac{(1 - r)^{2 + 2\alpha}(1-\alpha - 2r) + (1 + r)^{1 + 2\alpha}(1-\alpha + 3(1 + \alpha)r - 2r^2) + 2(1-\alpha + 2r)(1 - r^2)^{2+\alpha}\cos(\alpha\pi)}{((1 - r)^{2 + 2\alpha} + (1 + r)^{2 + 2\alpha} + 2(1 - r^2)^{1 + \alpha}\cos(\alpha\pi))^2}
\end{split}
\eqe
for $r\in(0,1)$ and $\Phi^Y_R(r)=0$ for $r>1$.
In Cartesian coordinates $H(\bol{x},t)$ can be calculated as
\eqb
H(\bol{x},t)=\frac{\Gamma(n+3/2)}{2\pi^{n+3/2} t\norm{\bol{x}}^{2n+2}}\Phi_R^Y\left(\frac{\norm{\bol{x}}}{t}\right),
\eqe
The density $H(\bol{x},t)$ is given by elementary functions for all $n$.
%
%
%
%
%
\subsection{Overshooting - odd number of dimensions $d=2n+3$}
For overshooting L\'evy walk we get \cite{MT}
\eqb \nonumber
H(\bol{k},s)=\frac{1}{s}g\left(\frac{i\bol{k}}{s}\right)=\frac{1}{s}\left(1-\frac{c\left(\frac{i\bol{k}}{s}\right)^\al(\cos(\pi\al/2)-i\sin(\pi\al/2))}{\int_{\bol{S}^{d-1}}\left(1-\scalar{\frac{i\bol{k}}{s}, \bol{u}}\right)^\alpha K(d\bol{u})}\right),
\eqe
where $c=\frac{\cos(\al \pi/2) \Gamma(2 - \al) \Gamma((1 + \al)/2) \Gamma(
  3/2 + n)}{(1 - \al) \sqrt{\pi} \Gamma(1 - \al) \Gamma(3/2 + \al/2 + n)}$
The inversion formula implies
\eqb
\begin{split}
&\Phi^Z_1(x)=\frac{c\Gamma(n+1)}{\Gamma(n+3/2)\sqrt{\pi} x^{2n+2+\al}}\mathcal{B}_1\cdot\bigg\{(1/x^2-1)^{n-m}\big[-\cos(\pi\al/2)\sin(\pi\al)(-1)^m(1/x-1)^{m+j+\al+1}\\
&+\sin(\pi\al/2)\left((-1)^j(1/x+1)^{m+j+\al+1}+(-1)^m(1/x-1)^{m+j+\al+1}\sin(\pi\al/2)\cos(\pi\al)\right)\big]\bigg\}\\
&/\mathcal{B}_2 \cdot\bigg\{\left((-1)^{j_1+j_2}(1+1/x)^{m_1+m_2+j_1+j_2+2\al+2}+(-1)^{m_1+m_2}(1/x-1)^{m_1+m_2+j_1+j_2+2\al+2} \right.\\
 &\left.+2\cos(\pi\al)(-1)^{j_1+m_2}(1/x+1)^{m_1+j_1+\al+1}(1/x-1)^{m_2+j_2+\al+1}\right)(1/x^2-1)^{2n-m_1-m_2}\bigg\}
\end{split}
\eqe
for $x\in(0,1)$ and
\eqb
\begin{split}
&\Phi^Z_1(x)=\frac{c\Gamma(n+1)\sin(\pi\al/2)}{\Gamma(n+3/2)\sqrt{\pi}x^{2n+2+\al}}
/\mathcal{B}_1 \cdot\bigg\{(-1)^{j}\left((1+1/x)^{m+j+\al+1}-(1-1/x)^{m+j+\al+1} \right)(1/x^2-1)^{n-m}\}
\end{split}
\eqe
for $x\in[1,\infty)$. For negative values of $x$ we take advantage from the symmetry $\Phi^Z_1(-x)=\Phi^Z_1(x)$. The symbols $\mathcal{B}_1$ and $\mathcal{B}_2$ are given by
\eqb
\mathcal{B}_1=\sum_{m}^n\sum_{j=0}^{m}{n \choose m}{m \choose j}2^{m-j}\frac{1}{m+j+\alpha+1},
\eqe
\eqb
\mathcal{B}_2=\sum_{m_1,m_2=0}^n\sum_{j_1=0}^{m_1}\sum_{j_2=0}^{m_2}{n \choose m_1}{m_1 \choose j_1}{n \choose m_2}{m_2 \choose j_2}2^{m_1+m_2-j_1-j_2}\frac{1}{(m_1+j_1+\alpha+1)(m_2+j_2+\alpha+1)}.
\eqe
Alternatively
\eqb
\begin{split}
\label{odd_over_phir_hyper}
\Phi^Z_1(x)=-\frac{c}{\pi \abs{x}^{\al+1}}\operatorname{Im}\frac{\cos(\pi\al/2)+i\sin(\pi\al/2)}{_2F_1(-\alpha/2, (1- \alpha)/2; 3/2+n; \frac{1}{x^2})}.
\end{split}
\eqe
The PDF $\Phi^Z_R$ is calculated from the equation
\eqb
\Phi_R^Z(\sqrt{r})=\frac{2\sqrt{\pi}}{\Gamma (n+3/2)}r^{n+1}(-1)^{n+1}\frac{d^{n+1}}{dr^{n+1}}\Phi_1^Z(\sqrt{r})
\eqe
if $r\in(0,\infty)$ and $\Phi^Z_R(r)=0$ otherwise.
In a special case $d=3$ we obtain
\eqb
\begin{split}
\Phi^Z_R(r)=&\frac{2\sin(\pi \al)}{\pi r}\Bigg ((r+1)^{2 \al+1} ((\al+2) r^2-3 (\al+1) r-1) (1-r)^\al\\
&\quad-2 (1-r^2)^\al (((\al+2) r+1) (1-r)^{\al+2}+(r+1)^{\al+2} ((\al+2) r-1)) \cos(\pi \al)\\
&\quad-(r+1)^\al (r (\al (r+3)+2 r+3)-1) (1-r)^{2 \al+1}\\
&\quad-(1-(\al+2) r) (1-r)^{3 \al+2}-(r+1)^{3 \al+2} ((\al+2) r+1))\Bigg)/ \\
&\quad\quad\quad\Bigg( (2 (1-r^2)^{\al+1} \cos(\pi \al)+(1-r)^{2 \al+2)}+(r+1)^{2 \al+2})^2\Bigg)
\end{split}
\eqe
for  $r\in (0,1)$ and
\eqb
\Phi^Z_R(r)=\frac{2\sin(\al\pi)}{\pi}\frac{ (1 + r)^\al(1 + (2 + \al)r)-(-1 + r)^\al(-1 + (2 + \al)r) }{
   r((-1 + r)^{1 + \al} - (1 + r)^{1 + \al})^2}
\eqe
for $r\in (1,\infty)$.
In Cartesian coordinates this gives us
\eqb
H(\bol{x},t)=\frac{\Gamma(n+3/2)}{2\pi^{n+3/2} t\norm{\bol{x}}^{2n+2}}\Phi_R^Z\left(\frac{\norm{\bol{x}}}{t}\right),
\eqe
The density $H(\bol{x},t)$ is given by elementary functions for all $n$.

\subsection{Undershooting - even number of dimensions $d=2n+2$}
In this case
\eqb
\begin{split}
\label{even_under_d_phir}
\Phi^Y_1(x)&=-\frac{1}{\pi \abs{x}}\operatorname{Im}\frac{1}{ _2F_1(-\alpha/2, (1- \alpha)/2; 1+n; \frac{1}{x^2})}
\end{split}
\eqe
on the interval $(0,1)$. Moreover
\eqb
\Phi_R^Y(\sqrt{r})=\frac{2\sqrt{\pi}}{\Gamma (n+1)}r^{n+1/2}D_-^{n+1/2}\{\Phi_1^Y(\sqrt{t})\}(r).
\eqe
if $r\in(0,1)$ and $\Phi^Y_R(r)=0$ otherwise.
Going back to Cartesian coordinates
\eqb
H(\bol{x},t)=\frac{\Gamma(n+1)}{2\pi^{n+1} t\norm{\bol{x}}^{2n+1}}\Phi_R^Y\left(\frac{\norm{\bol{x}}}{t}\right).
\eqe

\subsection{Overshooting - even number of dimensions $d=2n+2$}
Now we have
\eqb
\begin{split}
\label{even_over_d_phir}
\Phi^Z_1(x)&=-\frac{c}{\pi \abs{x}^{\al+1}}\operatorname{Im}\frac{\cos(\frac{1}{2}\pi\al)+i\sin(\frac{1}{2}\pi\al)}{ _2F_1(-\alpha/2, (1- \alpha)/2; 1+n; \frac{1}{x^2})}.
\end{split}
\eqe
In the above equation the constant $c$ is defined as $c=\frac{\cos(\al \pi/2) \Gamma(2 - \al) \Gamma((1 + \al)/2) \Gamma(
  1 + n)}{(1 - \al) \sqrt{\pi} \Gamma(1 - \al) \Gamma(1 + \al/2 + n)}$. The density of the radius:
\eqb
\Phi_R^Z(\sqrt{r})=\frac{2\sqrt{\pi}}{\Gamma (n+1)}r^{n+1/2}D_-^{n+1/2}\{\Phi_1^Z(\sqrt{t})\}(r).
\eqe
if $r\in(0,\infty)$ and $\Phi^Z_R(r)=0$ otherwise.
In Cartesian coordinates
\eqb
H(\bol{x},t)=\frac{\Gamma(n+1)}{2\pi^{n+1} t\norm{\bol{x}}^{2n+1}}\Phi_R^Z\left(\frac{\norm{\bol{x}}}{t}\right).
\eqe
\appendix

\section{Inversion of F-L transform for 1D ballistic processes}
In this appendix we present the method of inverting the Fourier-Laplace transform for 1D ballistic L\'evy walks from \cite{asymptotic densities} and \cite{godreche luck}. The Fourier-Laplace transform of $H_1(x,t)$ is given by
\begin{eqnarray}
H_1(k,s)&&=\int_{-\infty}^\infty\int_0^\infty e^{-ikx-st}H_1(x,t)dtdx\\
&&=\int_{-\infty}^\infty\int_0^\infty e^{-ikx-st}\frac{1}{t}\Phi_1\left(\frac{x}{t}\right)dtdx.
\end{eqnarray}
We substitute $y=\frac{x}{t}$ to get
\begin{eqnarray}
H_1(k,s)&&=\int_{-\infty}^\infty\int_0^\infty e^{-(iky+s)t}\Phi_1\left(y\right)dtdy
\\&&=\int_{-\infty}^\infty \frac{1}{iky+s}\Phi_1\left(y\right)dy=\frac{1}{s}\int_{-\infty}^\infty \frac{1}{\frac{ik}{s}y+1}\Phi_1\left(y\right)dy\\
&&=\frac{1}{s}\ex \frac{1}{\frac{ik}{s}X_1(1)+1}.
\end{eqnarray}
Hence
\eqb
\label{representation g_1}
g_1(\xi)=\ex \frac{1}{1+\xi X_1(1)}.
\eqe
Now, the density $\Phi_1(x)$ can be written as
\eqb
\label{phi ex}
\Phi_1(x)=\ex \delta(x-X_1(1)).
\eqe
The Sokhotsky-Weierstrass theorem implies
\eqb
\lim_{\eps\rightarrow 0}\frac{1}{x \pm i\eps}=\frac{1}{x}\mp i\pi \delta(x).
\eqe
We infer that
\eqb
\label{SW delta}
\delta(x)=-\operatorname{Im}\lim_{\eps\rightarrow 0}\frac{1}{x + i\eps}.
\eqe
Combining Eqs. \eqref{phi ex} and \eqref{SW delta} gives us
\begin{eqnarray}
\Phi_1(x)&&=-\operatorname{Im}\lim_{\eps\rightarrow 0}\ex\frac{1}{x-X_1(1) + i\eps}\\
&&=-\lim_{\eps\rightarrow 0}\operatorname{Im}\frac{1}{x+i\eps}\ex\frac{1}{1-\frac{X_1(1)}{x+i\eps}}
\end{eqnarray}
which yields
\eqb
\label{SW application}
\Phi_1(x)=-\frac{1}{\pi}\lim_{\eps \rightarrow 0} \operatorname{Im}\left[\frac{1}{x+i\eps}\ex \frac{1}{1-\frac{X_1(1)}{x+i\eps}}\right].
\eqe
Now we combine Eqs. \eqref{representation g_1} and \eqref{SW application} and get the formula for $\Phi_1$ (Eq. \eqref{formula_phi}).
 Notice that to formally get Eq. \eqref{formula_phi}  we have to choose the Fourier space variable $k$ and Laplace space variable $s$ such that $\frac{ik}{s}=\frac{-1}{x+i\eps}$. This is obtained by taking $s=s_1+is_2$ where $s_1=\eps>0$, $s_2=-x$ and $k=1$. Since $s_1>0$  and $k\in\R$ the Fourier-Laplace transform is well defined and $H_1(k,s)=\frac{1}{s}g_1\left(\frac{ik}{s}\right)$ even for complex values of $s$ (providing that $s_1>0$).

\section*{Acknowledgment}
This research was partially supported by NCN Maestro grant no. 2012/06/A/ST1/00258.

\end{document}